\documentclass[a4paper,11pt]{amsart}
\usepackage{graphicx}
\usepackage{mathptmx}
\usepackage{amsmath}
\usepackage{amssymb}
\usepackage{enumitem}
\usepackage{xcolor}
\usepackage{xparse}
\NewDocumentCommand{\eulerian}{omm}
 {%
  \genfrac<>{0pt}{}{#2}{#3}%
  \IfValueT{#1}{_{\!#1}}%
 }

 \newmuskip\pFqmuskip

\newcommand*\pFq[6][8]{%
  \begingroup % only local assignments
  \pFqmuskip=#1mu\relax
  \mathchardef\normalcomma=\mathcode`,
  \mathcode`\,=\string"8000
  \begingroup\lccode`\~=`\,
  \lowercase{\endgroup\let~}\pFqcomma
  {}_{#2}F_{#3}{\left(\genfrac..{0pt}{}{#4}{#5}\bigg|#6\right)}%
  \endgroup
}
\newcommand{\pFqcomma}{{\normalcomma}\mskip\pFqmuskip}

\newtheorem{theorem}{Theorem}
\newtheorem{lemma}[theorem]{Lemma}

\newtheorem{remark}[theorem]{Remark}

\begin{document}

\title[Some identities on degenerate Bell polynomials and their related identities]{Some identities on degenerate Bell polynomials and their related identities}

\author{Taekyun  Kim}
\address{Department of Mathematics, Kwangwoon University, Seoul 139-701, Republic of Korea}
\email{tkkim@kw.ac.kr}

\author{DAE SAN KIM}
\address{Department of Mathematics, Sogang University, Seoul 121-742, Republic of Korea}
\email{dskim@sogang.ac.kr}


\maketitle

\begin{abstract}
Recently, several types of degenerate Bell polynomials have been introduced as degenerate versions of the ordinary Bell polynomials. The aim of this paper is to study some identities for the degenerate Bell polynomials and their related identities by making use of series transformation formula obtained by Boyadzhiev. In particular, this answers to the natural question about the relationship between two different types of degenerate Bell polynomials.
\end{abstract}

\section{Introduction}
In [3], Carlitz initiated a study of degenerate Bernoulli and degenerate Euler polynomials, which are degenerate versions of the ordinary Bernoulli and Euler polynomials. In recent years, some mathematicians have explored degenerate versions of many special polynomials and numbers which include the degenerate Bernoulli numbers of the second kind, degenerate Stirling numbers of the first and second kinds, degenerate Bell numbers and polynomials, and so on (see [5--9,12] and the references therein). \par

The Bell number $B_n$ counts  the number of partitions of a set with $n$ elements into disjoint nonempty subsets. The Bell polynomials are natural extensions of Bell numbers and also called Touchard or exponential polynomials. As degenerate versions of these Bell polynomials and numbers, several different versions of degenerate Bell polynomials and numbers are introduced and investigated (see [7--9]).\par

The aim of this paper is to study some identities for the degenerate Bell polynomials and their related identities by using generating functions and the series transformation formula obtained by Boyadzhiev (see [2]). In particular, this answers to the natural question about the relationship between two different types of degenerate Bell polynomials. In addition, it gives us an expression for the generating function of the sums of $\lambda$-falling factorials of consecutive nonnegative integers, which reduces to the generating function of sums of powers of consecutive nonnegative integers. This paper demonstrates that the series transformation formula can be used in studying various types of the degenerate Bell polynomials as well.  \par

This paper is outlined as follows. In Section 1, we will recall degenerate exponential functions, degenerate Stirling numbers of the first kind and degenerate Stirling numbers of the second kind which are respectively degenerate versions of exponential functions, Stirling numbers of the first kind and Stirling numbers of the second kind. Then we will go over the definitions of some degenerate special polynomials, namely three types of degenerate Bell polynomials and the degenerate Bernoulli polynomials. As to two types degenerate Bell polynomials among the three types, we ask what their relationship is. In addition, we recall the geometric polynomials in connection with the ordinary Bell polynomials.

In Section 2, we will recall explicit expressions for the three types of degenerate Bell polynomials. Then we will introduce the degenerate geometric polynomilas in connection with the partially degenerate Bell polynomials. Next, we will state Theorem 3 about the series transformation formula which is proved in [2]. This is a key fact that we are going to apply in order to get some interesting results. By choosing $g(x)=e_{\lambda}(x)$ as the degenerate exponentials in Theorem 3, we get Theorem 4. And then, further taking $f(x)=e^{xt}$ as the usual exponential, we will be able to answer to the question regarding the relationship between two different types of degenerate Bell polynomials. As is well known, the sums of powers of consecutive nonnegative integers can be expressed in terms of Bernoulli polynomials. Here we express the generating function of the sums of $\lambda$-falling factorials of consecutive nonnegative integers in terms of the geometric polynomial and the Stirling numbers of the first kind. By letting $\lambda \rightarrow 0$, the generating function of the sums of powers of consecutive nonnegative integers is expressed by the geometric polynomial. Then we obtain an identity involving Stirling numbers of the second kind, the Carlitz degenerate Bernoulli numbers and the degenerate geometirc polynomial.
Finally, we conclude this paper in Section 3. \par

\vspace{0.1in}

For any nonzero $\lambda\in\mathbb{R}$, the degenerate exponential functions are defined by 
\begin{equation}
e_{\lambda}^{x}(t)=(1+\lambda t)^{\frac{x}{\lambda}},\quad e_{\lambda}(t)=e_{\lambda}^{1}(t)=(1+\lambda t)^{\frac{1}{\lambda}},\quad (\mathrm{see}\ [5]).\label{1}
\end{equation}
Thus, by \eqref{1}, we get 
\begin{equation}
e_{\lambda}^{x}(t)=\sum_{n=0}^{\infty}(x)_{n,\lambda}\frac{t^{n}}{n!},\quad(\mathrm{see}\ [6]),\label{2}
\end{equation}
where the $\lambda$-falling factorials are given by
\begin{equation}
(x)_{0,\lambda}=1,\quad (x)_{n,\lambda}=x(x-\lambda)(x-2\lambda)\cdots(x-(n-1)\lambda),\ n\ge 1. \label{3}
\end{equation}
It is well known that the Stirling numbers of the first kind are defined by 
\begin{equation}
(x)_{n}=\sum_{l=0}^{n}S_{1}(n,l)x^{l},\quad (n\ge 0),\quad(\mathrm{see}\ [7,9,11]),\label{4}
\end{equation}
where $(x)_{0}=1,\ (x)_{n}=x(x-1)\cdots(x-n+1),\ (n\ge 1)$. \par 
As the inversion of \eqref{4}, the Stirling numbers of the second kind are given by 
\begin{equation}
x^{n}=\sum_{l=0}^{n}S_{2}(n,l)(x)_{l},\quad(n\ge 0),\quad(\mathrm{see}\ [7,11]).\label{5}
\end{equation}
In [5], the degenerate Stirling numbers of the first kind are introduced as 
\begin{equation}
(x)_{n}=\sum_{l=0}^{n}S_{1,\lambda}(n,l)(x)_{n,\lambda},\quad (n\ge 0). \label{6}	
\end{equation}
As the inversion of \eqref{6}, the degenerate Stirling numbers of the second kind are given by 
\begin{equation}
(x)_{n,\lambda}=\sum_{l=0}^{n}S_{2,\lambda}(n,l)(x)_{l},\quad(n\ge 0),\quad(\mathrm{see}\ [5]).\label{7}
\end{equation}
Recently, the degenerate Bell polynomials are defined by 
\begin{equation}
e_{\lambda}\big(x(e^{t}-1)\big)=\sum_{n=0}^{\infty}\mathrm{Bel}_{n,\lambda}(x)\frac{t^{n}}{n!},\quad(\mathrm{see}\ [9]).\label{8}
\end{equation}
Note that $\displaystyle\lim_{\lambda\rightarrow 0}\mathrm{Bel}_{n,\lambda}(x)=\mathrm{Bel}_{n}(x),\ (n\ge 0)\displaystyle$, where $\mathrm{Bel}_{n}(x)$ are the ordinary Bell polynomials given by 
\begin{equation}
e^{x(e^{t}-1)}=\sum_{n=0}^{\infty}\mathrm{Bel}_{n}(x)\frac{t^{n}}{n!},\quad(\mathrm{see}\ [1,11]).\label{9}
\end{equation}
In [7], the degenerate Bell polynomials of the second kind are defined by 
\begin{equation}
e_{\lambda}(xe^{t})\cdot e_{\lambda}^{-1}(x)=\sum_{n=0}^{\infty}\mathrm{bel}_{n,\lambda}(x)\frac{t^{n}}{n!}.\label{10}	
\end{equation}
Observe that $\displaystyle\lim_{\lambda\rightarrow 0}\mathrm{bel}_{n,\lambda}(x)=\mathrm{Bel}_{n}(x) \ (n\ge 0)\displaystyle$.
The following natural question arises natually from \eqref{8} and \eqref{10}.\\
\\~~
\textbf{Question.} What is the relationship between $\mathrm{Bel}_{n,\lambda}(x)$ and $\mathrm{bel}_{n,\lambda}(x)$? \\
\\~~
This question will be answered in Theorem 6. \par 

For $n\ge 0$, the geometric polynomials are defined by 
\begin{equation}
W_{n}(x)=\sum_{k=0}^{n}S_{2}(n,k)k!x^{k},\quad(\mathrm{see}\ [2,12]).\label{11}	
\end{equation}
From \eqref{11}, we note that 
\begin{equation}
\sum_{n=0}^{\infty}W_{n}(x)\frac{t^{n}}{n!}=\frac{1}{1-x(e^{t}-1)}.\label{12}
\end{equation}
By \eqref{12}, we get 
\begin{equation}
W_{n}(x)=\int_{0}^{\infty}\mathrm{Bel}_{n}(xy)e^{-y}dy,\quad(\mathrm{see}\ [2]).\label{13}
\end{equation}
In [8], the partially degenerate Bell polynomials are defined by 
\begin{equation}
e^{x(e_{\lambda}(t)-1)}=\sum_{n=0}^{\infty}\phi_{n,\lambda}(x)\frac{t^{n}}{n!}.\label{14}	
\end{equation}
Note that $\displaystyle\lim_{\lambda\rightarrow 0}\phi_{n,\lambda}(x)=\mathrm{Bel}_{n}(x),\ (n\ge 0)\displaystyle$. 
\par 
Carlitz introduced the degenerate Bernoulli numbers given by 
\begin{equation}
\frac{t}{e_{\lambda}(t)-1}=\sum_{n=0}^{\infty}\beta_{n,\lambda}\frac{t^{n}}{n!},\quad(\mathrm{see}\ [3]).\label{15}
\end{equation}
Note that $\displaystyle\lim_{\lambda\rightarrow 0}\beta_{n,\lambda}=B_{n},\ (n\ge 0)\displaystyle$, where $B_{n}$ are the ordinary Bernoulli numbers given by
\begin{equation}
\frac{t}{e^{t}-1}=\sum_{n=0}^{\infty}B_{n}\frac{t^{n}}{n!},\quad(\mathrm{see}\ [1-12]).\label{16}
\end{equation}
More generally, the Bernoulli polynomials $B_{n}(x)$ are defined by
\begin{equation}
\frac{t}{e^{t}-1}e^{xt}=\sum_{n=0}^{\infty}B_{n}(x)\frac{t^{n}}{n!}. \label{17}
\end{equation}

\section{Some identities on degenerate Bell polynomials and their related identities}
From \eqref{8} and \eqref{10}, it is immediate to see that 
\begin{equation}
\mathrm{Bel}_{n,\lambda}(x)=\sum_{k=0}^{n}(1)_{k,\lambda}S_{2}(n,k)x^{k},\quad(n\ge 0),\label{18}
\end{equation}
and 
\begin{equation}
\mathrm{bel}_{n,\lambda}(x)=e_{\lambda}^{-1}(x)\sum_{k=0}^{\infty}\frac{(1)_{k,\lambda}}{k!}k^{n}x^{k},\quad(n\ge 0).\label{19}	
\end{equation}
By \eqref{14}, we get 
\begin{equation}
\phi_{n,\lambda}(x)=\sum_{k=0}^{n}S_{2,\lambda}(n,k)x^{k},\quad(n\ge 0).\label{20}	
\end{equation}
By inversion, from \eqref{18} and \eqref{20}, we have
\begin{align*}
x^{n}=\frac{1}{(1)_{n,\lambda}}\sum_{k=0}^{n}S_{1}(n,k)\mathrm{Bel}_{k,\lambda}(x), \quad
x^{n}=\sum_{k=0}^{n}S_{1,\lambda}(n,k)\phi_{k,\lambda}(x).
\end{align*}
Now, we define the degenerate geometric polynomials by 
\begin{equation}
W_{n,\lambda}(x)=\sum_{k=0}^{n}S_{2,\lambda}(n,k)k!	x^{k},\quad(n\ge 0).\label{21}
\end{equation}
Note that $\displaystyle\lim_{\lambda\rightarrow 0}W_{n,\lambda}(x)=W_{n}(x),\ (n\ge 0)\displaystyle$. Further, by inversion we also have
\begin{equation*}
x^{n}=\frac{1}{n!}\sum_{k=0}^{n}S_{1,\lambda}(n,k)W_{k,\lambda}(x).
\end{equation*}
From \eqref{21}, we have 
\begin{align}
\sum_{n=0}^{\infty}W_{n,\lambda}(x)\frac{t^{n}}{n!}\ &=\ \sum_{n=0}^{\infty}\bigg(\sum_{k=0}^{n}k!S_{2,\lambda}(n,k)x^{k}\bigg)\frac{t^{n}}{n!} \label{22}\\
&=\ \sum_{k=0}^{\infty}	x^{k}k!\sum_{n=k}^{\infty}S_{2,\lambda}(n,k)\frac{t^{n}}{n!}\nonumber \\
&=\ \sum_{k=0}^{\infty}x^{k}\big(e_{\lambda}(t)-1\big)^{k}=\frac{1}{1-x(e_{\lambda}(t)-1)}. \nonumber
\end{align}
We observe that 
\begin{align}
\int_{0}^{\infty}e^{xy(e_{\lambda}(t)-1)}e^{-y}dy&=\int_{0}^{\infty}e^{-y(1-x(e_{\lambda}(t)-1))}dy \label{23} \\
&=\frac{1}{1-x(e_{\lambda}(t)-1)}=\sum_{n=0}^{\infty}W_{n,\lambda}(x)\frac{t^{n}}{n!}.\nonumber
\end{align}
On the other hand,
\begin{equation}
\int_{0}^{\infty}e^{xy(e_{\lambda}(t)-1)}e^{-y}dy=\sum_{n=0}^{\infty}\int_{0}^{\infty}\phi_{n,\lambda}(xy)e^{-y}dy\frac{t^{n}}{n!}.\label{24}
\end{equation}
Therefore, by \eqref{23} and \eqref{24}, we obtain the following theorem. 
\begin{theorem}
For $n\ge 0$, we have 
\begin{displaymath}
\int_{0}^{\infty}\phi_{n,\lambda}(xy)e^{-y}dy=W_{n,\lambda}(x). 
\end{displaymath}	
\end{theorem}
Later, we will need the following lemma which is proved in [4, p.218] and also mentioned in [2].
\begin{lemma}
Let $D=\frac{d}{dx}$, and let $f$ be analytic on an open set $U$ in $\mathbb{C}$. Then, for each $x \in U$, we have 
\begin{equation*}
(xD)^{n}f(x)=\sum_{k=0}^{n}S_{2}(n,k)x^{k}D^{k}f(x),\quad(n\ge 0).
\end{equation*}
\end{lemma}
In [2], K. N. Boyadzhiev presented a formula that turns power series into series of functions. He also gave many valuable formulas related to polynomials. 
The following series transformation formula is proved in [2, Theorem 4.1]. Indeed, it can be shown by computing $f(xD)g(x)$ in two different ways with the help of Lemma 2. We are going to apply this fromula for various $f$ and $g$ in order to get a few interesting results.
\begin{theorem}
Assume that $f(x)=\sum_{n=0}^{\infty}a_{n}x^{n}$ and  $g(x)=\sum_{k=0}^{\infty}c_{k}x^{k}$ are power series convergnt on some open disks centered at the origin. Then we have 
\begin{align*}
\sum_{k=0}^{\infty}\frac{g^{(k)}(0)}{k!}f(k)x^{k}=\sum_{n=0}^{\infty}\frac{f^{(n)}(0)}{n!}\sum_{k=0}^{n}S_{2}(n,k)g^{(k)}(x)x^{k}.
\end{align*}
\end{theorem}

First, we choose $g$ as $g(x)=e_{\lambda}(x)$ in Theorem 3. Then we see that the following holds:\par 
\begin{equation}
g^{(k)}(x)=\bigg(\frac{d}{dx}\bigg)^{k}e_{\lambda}(x)=(1)_{k,\lambda}e_{\lambda}^{1-k \lambda}(x)=\frac{(1)_{k,\lambda}}{(1+\lambda x)^{k}}e_{\lambda}(x).\label{25}	
\end{equation}
From Theorem 3 and \eqref{25}, we note that 
\begin{equation}
\sum_{n=0}^{\infty}\frac{(1)_{n,\lambda}}{n!}f(n)x^{n}=e_{\lambda}(x)\sum_{n=0}^{\infty}\frac{f^{(n)}(0)}{n!}\sum_{k=0}^{n}S_{2}(n,k)(1)_{k,\lambda}\bigg(\frac{x}{1+\lambda x}\bigg)^{k}.\label{26}	
\end{equation}
By \eqref{18} and \eqref{26}, we get 
\begin{equation}
\sum_{n=0}^{\infty}\frac{(1)_{n,\lambda}}{n!}f(n)x^{n}=e_{\lambda}(x)\sum_{n=0}^{\infty}\frac{f^{(n)}(0)}{n!}\mathrm{Bel}_{n,\lambda}\bigg(\frac{x}{1+\lambda x}\bigg).\label{27}
\end{equation}
Therefore, we obtain the following theorem. 
\begin{theorem}
Let $f(x)=\sum_{n=0}^{\infty}a_{n}x^{n}$ be a power series convergnt on some open disk centered at the origin. For $n\ge 0$, we have 
\begin{displaymath}
\sum_{n=0}^{\infty}\frac{(1)_{n,\lambda}}{n!}f(n)x^{n}=e_{\lambda}(x)\sum_{n=0}^{\infty}\frac{f^{(n)}(0)}{n!}\mathrm{Bel}_{n,\lambda}\bigg(\frac{x}{1+\lambda x}\bigg).
\end{displaymath}
\end{theorem}
Now, by using \eqref{19} we observe that 
\begin{align}
(xD)^{n}e_{\lambda}(x) &=(xD)^{n}\sum_{k=0}^{\infty}\frac{(1)_{k,\lambda}}{k!}x^{k}=\sum_{k=0}^{\infty}\frac{(1)_{k,\lambda}}{k!}k^{n}x^{k}\label{28}\\
&=\bigg(e_{\lambda}^{-1}(x)\sum_{k=0}^{\infty}\frac{k^{n}}{k!}(1)_{k,\lambda}x^{k}\bigg)e_{\lambda}(x)=\mathrm{bel}_{n,\lambda}(x)e_{\lambda}(x). \nonumber
\end{align}
Hence, by \eqref{28}, we get the next result.
\begin{theorem}
For $n\ge 0$, we have 	
\begin{displaymath}
\mathrm{bel}_{n,\lambda}(x)e_{\lambda}(x)=(xD)^{n}e_{\lambda}(x)=\sum_{k=0}^{\infty}\frac{(1)_{k,\lambda}}{k!}k^{n}x^{k}.
\end{displaymath}
\end{theorem}
Taking $f(x)=e^{xt}$ in \eqref{26}, we have 
\begin{align}
\sum_{n=0}^{\infty}\frac{(1)_{n,\lambda}}{n!}e^{nt}x^{n}&=e_{\lambda}(x)\sum_{n=0}^{\infty}\frac{t^{n}}{n!}\sum_{k=0}^{n}S_{2}(n,k)(1)_{k,\lambda}\bigg(\frac{x}{1+\lambda x}\bigg)^{k}\label{29} \\
&=e_{\lambda}(x)\sum_{n=0}^{\infty}\mathrm{Bel}_{n,\lambda}\bigg(\frac{x}{1+\lambda x}\bigg)\frac{t^{n}}{n!}.\nonumber
\end{align}
Thus, by \eqref{29}, we get 
\begin{align}
\sum_{n=0}^{\infty}\mathrm{Bel}_{n,\lambda}\bigg(\frac{x}{1+\lambda x}\bigg)\frac{t^{n}}{n!}&=e_{\lambda}^{-1}(x)\sum_{n=0}^{\infty}\frac{(1)_{n,\lambda}}{n!}e^{nt}x^{n}\label{30} \\
&=e_{\lambda}^{-1}(x)e_{\lambda}\big(xe^{t}\big)=\sum_{n=0}^{\infty}\mathrm{bel}_{n,\lambda}(x)\frac{t^{n}}{n!}.\nonumber	
\end{align}
Therefore, by comparing the coefficients on both sides of \eqref{30}, we obtain the following theorem. This answers to the question that is asked in the Introduction.
\begin{theorem}
For $n\ge 0$, we have 
\begin{displaymath}
\mathrm{bel}_{n,\lambda}(x)=\mathrm{Bel}_{n,\lambda}\bigg(\frac{x}{1+\lambda x}\bigg).
\end{displaymath}	
\end{theorem}
From \eqref{8}, \eqref{10} and \eqref{30}, we note that 
\begin{equation}
e_{\lambda}\bigg(\frac{x}{1+\lambda x}\big(e^{t}-1\big)\bigg)=e_{\lambda}(xe^{t})\cdot e_{\lambda}^{-1}(x).\label{31}
\end{equation}
For $g(x)=\frac{1}{1-x}$, we have $g^{(k)}(x)=\big(\frac{d}{dx}\big)^{k}g(x)=k!\big(\frac{1}{1-x}\big)^{k+1}$, $g^{(k)}(0)=k!$. \par 
From Theorem 3 , we have 
\begin{align}
\sum_{k=0}^{\infty}\frac{k!}{k!}f(k)x^{k}&=\sum_{n=0}^{\infty}\frac{f^{(n)}(0)}{n!}\sum_{k=0}^{n}S_{2}(n,k)x^{k}k!\bigg(\frac{1}{1-x}\bigg)^{k+1}\label{32} \\
&=\frac{1}{1-x}\sum_{n=0}^{\infty}\frac{f^{(n)}(0)}{n!}\sum_{k=0}^{n}S_{2}(n,k)\bigg(\frac{x}{1-x}\bigg)^{k}k!
\nonumber\\
&=\frac{1}{1-x} \sum_{n=0}^{\infty}\frac{f^{(n)}(0)}{n!}W_{n}\bigg(\frac{x}{1-x}\bigg).\nonumber
\end{align}
Thus, by \eqref{32}, we get 
\begin{equation}
\sum_{k=0}^{\infty}f(k)x^{k}=\frac{1}{1-x}\sum_{n=0}^{\infty}\frac{f^{(n)}(0)}{n!}W_{n}\bigg(\frac{x}{1-x}\bigg).\label{33}	
\end{equation}
Applying \eqref{33} to $f(x)=(x)_{m,\lambda}$, and noting 
\begin{displaymath}
f(x)=(x)_{m,\lambda}=\lambda^{m}\bigg(\frac{x}{\lambda}\bigg)_{m}=\sum_{l=0}^{m}S_{1}(m,l)\lambda^{m-l}x^{l},
\end{displaymath}
and, for $ n \le m$,
\begin{equation}
f^{(n)}(x)=\bigg(\frac{d}{dx}\bigg)^{n}f(x)=\sum_{l=n}^{m}S_{1}(m,l)\lambda^{m-l}(l)_{n}x^{l-n}.\label{34}
\end{equation}
From \eqref{33} and \eqref{34}, we have 
\begin{equation}
\sum_{k=0}^{\infty}(k)_{m,\lambda}x^{k}=\frac{1}{1-x}\sum_{l=0}^{m}S_{1}(m,l)\lambda^{m-l}W_{l}\bigg(\frac{x}{1-x}\bigg).\label{35}	
\end{equation}
We observe that 
\begin{align}
x^{m\lambda}\big(x^{1-\lambda}D\big)^{m}\bigg(\frac{1}{1-x}\bigg)&=x^{m\lambda}\big(x^{1-\lambda}D\big)^{m}\sum_{k=0}^{m}x^{k} \label{36} \\
&= \sum_{k=0}^{\infty}(k)_{m,\lambda}x^{k}.  \nonumber
\end{align}
From \eqref{35} and \eqref{36}, we have 
\begin{align}
x^{m\lambda}\big(x^{1-\lambda}D\big)^{m}\bigg(\frac{1}{1-x}\bigg)&=\sum_{k=0}^{\infty}(k)_{m,\lambda}x^{k}\label{37} \\
&=\frac{1}{1-x}\sum_{l=0}^{m}S_{1}(m,l)\lambda^{m-l}W_{l}\bigg(\frac{x}{1-x}\bigg). \nonumber
\end{align}
Note that 
\begin{align}
\frac{1}{1-x}\sum_{l=0}^{\infty}(l)_{m,\lambda}x^{l}&=\bigg(\sum_{j=0}^{\infty}x^{j}\bigg)\bigg(\sum_{l=0}^{\infty}(l)_{m,\lambda}x^{l}\bigg)=\sum_{k=0}^{\infty}\bigg(\sum_{l=0}^{k}(l)_{m,\lambda}\bigg)x^{k}\label{38} \\
&=\sum_{k=0}^{\infty}\big((0)_{m,\lambda}+(1)_{m,\lambda}+(2)_{m,\lambda}+\cdots+(k)_{m,\lambda}\big)x^{k}.\nonumber 
\end{align}
By \eqref{38}, we get 
\begin{align}
\sum_{k=0}^{\infty}\big((0)_{m,\lambda}+(1)_{m,\lambda}+(2)_{m,\lambda}+\cdots+(k)_{m,\lambda}\big)x^{k}&=\frac{1}{1-x}\sum_{l=0}^{\infty}(l)_{m,\lambda}x^{l} \label{39} \\
&=\frac{1}{(1-x)^{2}}\sum_{l=0}^{m}S_{1}(m,l)\lambda^{m-l}W_{l}\bigg(\frac{x}{1-x}\bigg).\nonumber
\end{align}
Therefore, by \eqref{39}, we obtain the following theorem. 
\begin{theorem}
For $m\ge 0$, we have 
\begin{displaymath}
\sum_{k=0}^{\infty}\big((0)_{m,\lambda}+(1)_{m,\lambda}+(2)_{m,\lambda}+\cdots+(k)_{m,\lambda}\big)x^{k}= \frac{1}{(1-x)^{2}}\sum_{l=0}^{m}S_{1}(m,l)\lambda^{m-l}W_{l}\bigg(\frac{x}{1-x}\bigg).
\end{displaymath}	
\end{theorem}
\begin{remark}
\textnormal{(a)} The above theorem is equivalent to the following:
\begin{align*}
&\sum_{k=1}^{\infty}\big((1)_{m,\lambda}+(2)_{m,\lambda}+\cdots+(k)_{m,\lambda}\big)x^{k}= \frac{1}{(1-x)^{2}}\sum_{l=0}^{m}S_{1}(m,l)\lambda^{m-l}W_{l}\bigg(\frac{x}{1-x}\bigg),\quad(m \ge 1), \\
&\quad\quad\quad\quad\quad \sum_{k=0}^{\infty}(k+1)x^{k}=\frac{1}{(1-x)^{2}}.
\end{align*}
\textnormal{(b)} It is well known that Faulhaber's formula on the sums of powers of consecutive nonnegative integers is equivalent to 
\begin{align}
\sum_{i=0}^{k} i^{m}=\frac{1}{m+1}\big(B_{m+1}(k+1)-B_{m+1}\big).\label{40}
\end{align}
By letting $\lambda \rightarrow 0$ in Theorem 7 and making use of \eqref{40}, we have
\begin{align*}
\sum_{k=0}^{\infty}\big(0^{m}+1^{m}+2^{m}+\cdots k^{m}\big)x^{k}&=\frac{1}{m+1}\sum_{k=0}^{\infty}\big(B_{m+1}(k+1)-B_{m+1}\big)x^{k} \\
&=\frac{1}{(1-x)^{2}}W_{m}\bigg(\frac{x}{1-x}\bigg).
\end{align*}
\end{remark}

For $n,m\in\mathbb{Z}$, with $ n \ge 1, 0 \le m \le n$, let $\eulerian{n}{m}$ be the Eulerian numbers. Then we note that 
\begin{equation}
\eulerian{n}{k}	=\sum_{l=0}^{k+1}(-1)^{l}\binom{k+1}{l}(k+1-l)^{n}\quad\mathrm{and}\quad \sum_{k=0}^{n}\eulerian{n}{k}=n!.\label{41}
\end{equation}
The Eulerian polynomials are defined by 
\begin{equation}
A_{n}(x)=\sum_{m=0}^{n}\eulerian{n}{m}x^{n-m},\quad (n\ge 0). \label{42}	
\end{equation}
From \eqref{37}, we note that 
\begin{equation}
\sum_{k=0}^{\infty}k^{m}x^{k}=\lim_{\lambda\rightarrow 0}\sum_{k=0}^{\infty}(k)_{m,\lambda}x^{k}=\frac{1}{1-x}W_{m}\bigg(\frac{x}{1-x}\bigg).\label{43}
\end{equation}
It is known that 
\begin{equation}
\sum_{k=0}^{\infty}k^{m}x^{k}=\frac{1}{(1-x)^{m+1}}A_{m}(x),\quad (\mathrm{see}\ [5]).\label{44}
\end{equation}
By \eqref{43} and \eqref{44}, we get 
\begin{equation}
\frac{1}{1-x}W_{m}\bigg(\frac{x}{1-x}\bigg)=\frac{1}{(1-x)^{m+1}}A_{m}(x).\label{45}	
\end{equation}
From \eqref{45}, we note that 
\begin{equation}
A_{m}(x)=(1-x)^{m}W_{m}\bigg(\frac{x}{1-x}\bigg).\label{46}
\end{equation}
By \eqref{37} and \eqref{46}, we get 
\begin{align}
\sum_{k=0}^{\infty}(k)_{m,\lambda}x^{k}&=\frac{1}{1-x}\sum_{l=0}^{m}S_{1}(m,l)\lambda^{m-l}W_{l}\bigg(\frac{x}{1-x}\bigg)\label{47} \\
&=\sum_{l=0}^{m}S_{1}(m,l)\lambda^{m-l}\frac{1}{(1-x)^{l+1}}A_{l}(x).	\nonumber
\end{align}
For $g(x)=(1-x)^{-r},\ (r\in\mathbb{N})$, we note that 
\begin{equation}
g^{(k)}(x)=\bigg(\frac{d}{dx}\bigg)^{k}g(x)=r(r+1)\cdots(r+k-1)(1-x)^{-r-k},\label{48}
\end{equation}
and hence that
\begin{align}
g^{(k)}(0)=r(r+1)\cdots(r+k-1)=\binom{-r}{k}k!(-1)^{k}.\label{49}
\end{align}
By Theorem 3, \eqref{48} and \eqref{49}, we get 
\begin{align}
\sum_{k=0}^{\infty}\frac{(-1)^{k}}{k!}k!\binom{-r}{k}f(k)x^{k}&=\sum_{n=0}^{\infty}\frac{f^{(n)}(0)}{n!}\sum_{k=0}^{n}S_{2}(n,k)r(r+1)\cdots(r+k-1)\frac{x^{k}}{(1-x)^{r+k}}\label{50}\\
&=\frac{1}{(1-x)^{r}}\sum_{n=0}^{\infty}\frac{f^{(n)}(0)}{n!}\sum_{k=0}^{n}S_{2}(n,k)\frac{\Gamma(r+k)}{\Gamma(r)}\bigg(\frac{x}{1-x}\bigg)^{k}\nonumber \\
&=\frac{1}{(1-x)^{r}}\sum_{n=0}^{\infty}\frac{f^{(n)}(0)}{n!}W_{n}\bigg(\frac{x}{1-x}\ \bigg|\ r\bigg), \nonumber
\end{align}
where 
\begin{equation}
W_{n}(x|r)=\sum_{k=0}^{n}S_{2}(n,k)\frac{\Gamma(r+k)}{\Gamma(k)}x^{k}.\label{51}	
\end{equation}
For $f(x)=(x)_{m,\lambda}$, we note that, for $n \le m$, 
\begin{equation}
\begin{aligned}
	\frac{1}{n!}\bigg(\frac{d}{dx}\bigg)^{n}f(x)&=\frac{1}{n!}f^{(n)}(x)=\frac{\lambda^{m}}{n!}\bigg(\frac{d}{dx}\bigg)^{n}\bigg(\frac{x}{\lambda}\bigg)_{m}\\
	&=\frac{1}{n!}\sum_{l=n}^{m}S_{1}(m,l)\lambda^{m-l}(l)_{n}x^{l-n}.
\end{aligned}	\label{52}
\end{equation}
By \eqref{50} and \eqref{52}, we get 
\begin{align}
\sum_{k=0}^{\infty}\binom{-r}{k}(-1)^{k}(k)_{m,\lambda}x^{k}&=\frac{1}{(1-x)^{r}}\sum_{n=0}^{m}\frac{1}{n!}\bigg(\sum_{l=0}^{m}S_{1}(m,l)\lambda^{m-l}(l)_{n}\bigg)0^{l-n}W_{n}\bigg(\frac{x}{1-x}\bigg|r\bigg)\label{53} \\
&=\frac{1}{(1-x)^{r}}\sum_{l=0}^{m}S_{1}(m,l)\lambda^{m-l}W_{l}\bigg(\frac{x}{1-x}\bigg|r\bigg).\nonumber
\end{align}
From \eqref{53}, we note that 
\begin{equation}
\sum_{k=0}^{\infty}\binom{-r}{k}(-1)^{k}(k)_{m,\lambda}x^{k}=\frac{1}{(1-x)^{r}}\sum_{l=0}^{m}S_{1}(m,l)\lambda^{m-l}W_{l}\bigg(\frac{1}{1-x}\bigg|r\bigg). \label{54}	
\end{equation}
On the other hand, 
\begin{align}
x^{\lambda m}\big(x^{1-\lambda}D\big)^{m}\frac{1}{(1-x)^{r}}&=x^{\lambda m}\big(x^{1-\lambda}D\big)^{m}\sum_{k=0}^{\infty}\binom{-r}{k}(-1)^{k}x^{k} \label{55} \\
&=\sum_{k=0}^{\infty}\binom{-r}{k}(-1)^{k}(k)_{m,\lambda}x^{k}.\nonumber
\end{align}
Thus, by \eqref{54} and \eqref{55}, we get 
\begin{equation}
x^{\lambda m}\big(x^{1-\lambda}D\big)^{m}\bigg(\frac{1}{1-x}\bigg)^{r}=\frac{1}{(1-x)^{r}}\sum_{l=0}^{m}	S_{1}(m,l)\lambda^{m-l}W_{l}\bigg(\frac{x}{1-x}\bigg|r\bigg).\label{56}
\end{equation}
Now, we observe that 
\begin{align}
\frac{1}{e_{\lambda}(t)+1}&=\frac{1}{e_{\lambda}(t)-1}-\frac{2}{e_{\lambda}^{2}(t)-1}\label{57}	\\
&=\frac{1}{t}\bigg(\frac{t}{e_{\lambda}(t)-1}-\frac{2t}{e_{\lambda}^{2}(t)-1}\bigg)=\frac{1}{t}\bigg(\frac{t}{e_{\lambda}(t)-1}-\frac{2t}{e_{\frac{\lambda}{2}}(2t)-1}\bigg)\nonumber \\
&=\frac{1}{t}\sum_{n=0}^{\infty}\big(\beta_{n,\lambda}-2^{n}\beta_{n,\frac{\lambda}{2}}\big)\frac{t^{n}}{n!}=\sum_{n=0}^{\infty}\bigg(\frac{\beta_{n+1,\lambda}-2^{n+1}\beta_{n+1,\frac{\lambda}{2}}}{n+1}\bigg)\frac{t^{n}}{n!}. \nonumber
\end{align}
Taking $x=-\frac{1}{2}$ in \eqref{21}, we have 
\begin{equation}
\sum_{n=0}^{\infty}W_{n,\lambda}\bigg(-\frac{1}{2}\bigg)\frac{t^{n}}{n!}=\frac{2}{e_{\lambda}(t)+1}=\sum_{n=0}^{\infty}2\bigg(\frac{\beta_{n+1,\lambda}-2^{n+1}\beta_{n+1,\frac{\lambda}{2}}}{n+1}\bigg)\frac{t^{n}}{n!}.\label{58}
\end{equation}
By comparing the coefficients on both sides of \eqref{58}, we get 
\begin{equation}
W_{n,\lambda}\bigg(-\frac{1}{2}\bigg)= 2\bigg(\frac{\beta_{n+1,\lambda}-2^{n+1}\beta_{n+1,\frac{\lambda}{2}}}{n+1}\bigg).\label{59}	
\end{equation}
From \eqref{20}, we have 
\begin{equation}
W_{n,\lambda}\bigg(-\frac{1}{2}\bigg)=\sum_{k=0}^{n}S_{2,\lambda}(n,k)k!(-1)^{k}2^{-k}\label{60}.
\end{equation}
Therefore, by \eqref{59} and \eqref{60}, we obtain the following theorem. 
\begin{theorem}
For $n\ge 0$, we have 
\begin{displaymath}
\sum_{k=0}^{n}S_{2,\lambda}(n,k)2^{-k}(-1)^{k}k!=W_{n,\lambda}\bigg(-\frac{1}{2}\bigg)=\frac{2}{n+1}\big(\beta_{n+1,\lambda}-2^{n+1}\beta_{n+1,\frac{\lambda}{2}}\big).
\end{displaymath}	
\end{theorem}
\begin{remark} 
We observe that 
\begin{align}
\sum_{n=0}^{\infty}\bigg(\frac{\beta_{n+1,\lambda}-2^{n+1}\beta_{n+1,\frac{\lambda}{2}}}{n+1}\bigg)\frac{t^{n}}{n!}&=\frac{1}{e_{\lambda}(t)+1}=\sum_{k=0}^{\infty}(-1)^{k}e_{\lambda}^{k}(t)\label{61} \\
&=\sum_{n=0}^{\infty}\bigg(\sum_{k=0}^{\infty}(-1)^{k}(k)_{n,\lambda}\bigg)\frac{t^{n}}{n!}.\nonumber 
\end{align}
By comparing the coefficients on both sides of \eqref{61}, we get 
\begin{displaymath}
\sum_{k=0}^{\infty}(-1)^{k}(k)_{n,\lambda}=\frac{1}{n+1}\big(\beta_{n+1,\lambda}-2^{n+1}\beta_{n+1,\frac{\lambda}{2}}\big). 
\end{displaymath}
\end{remark}

\section{Conclusion}

It was Carlitz who initiated a study of the degenerate Bernoulli and degenerate Euler polynomials. In recent years, intensive studies have been done for degenerate versions of quite a few special polynomials and numbers. Here we investigated three types of Bell polynomials by using the series transformation formula proved by Boyadzhiev in [2]. We obtained some interesting results by applying this key result to various choices of power series $f$ and $g$. Specifically, we answered to the natural question regarding the relationship between two different types of degenerate Bell polynomials. In addition, we expressed the generating function of the sums of $\lambda$-falling factorials of consecutive nonnegative integers in terms of the geometric polynomial and the Stirling numbers of the first kind. \par
It is one of our future research projects to continue this line of research, namely to study various degenerate versions of many special polynomials and numbers and their applications to science, physics and engineering.

\end{document}